\documentclass[amssymb,twocolumn,aps,nofootinbib]{revtex4}
\usepackage{amsmath}
\usepackage{amstext}
\usepackage{amsopn}
\usepackage{amsfonts}
\usepackage{amssymb}
\usepackage{amssymb}
\usepackage{bbm}
\usepackage{accents}
\usepackage{empheq}
\usepackage{graphicx}
\usepackage{epsf}
\usepackage{graphics}
\usepackage{xcolor}
\usepackage[latin1]{inputenc}
\def\sc{\scriptstyle}

\begin{document}

\title{Supersymmetry and the Riemann zeros on the critical line}

\author{Ashok Das$^{a}$ and Pushpa Kalauni$^{b}$}
\affiliation{$^a$ Department of Physics and Astronomy, University of Rochester,
Rochester, NY 14627-0171, USA}
\affiliation{$^b$ Homer L. Dodge Department of Physics and Astronomy, University
of Oklahoma, Norman, OK 73019, USA}

\begin{abstract}

We propose a new way of studying the Riemann zeros on the critical line using ideas from supersymmetry. Namely, we construct a supersymmetric quantum mechanical model whose energy eigenvalues correspond to the Riemann zeta function in the strip $0< {\rm Re}\, s < 1$ (in the complex parameter space) and show that the zeros on the critical line arise naturally from the vanishing ground state energy condition in this model. 

\end{abstract}
\maketitle

Riemann \cite{riemann} generalized Euler's zeta function to the entire complex space of parameters in three essential  steps. First, he extended the series representation of the zeta function to complex parameters as
\begin{equation}
\zeta(s) = \sum_{n=1}^{\infty} \frac{1}{n^{s}},\ \ s=\sigma + i\lambda,\ (\sigma,\lambda\,{\rm real}),\quad {\rm Re}\,s>1.\label{1}
\end{equation}
This series representation can be written as a product of factors involving only prime numbers and the zeta function in \eqref{1}  has an integral representation in terms of the Mellin transform of the Bose-Einstein distribution function \cite{abramowitz}. In the second step, he expressed the zeta function in \eqref{1} in terms of the alternating zeta function as
\begin{align}
\zeta(s) & = \frac{1}{1-2^{1-s}}\, \eta(s)\notag\\
& = \frac{1}{1-2^{1-s}} \sum_{n=1}^{\infty} \frac{(-1)^{n+1}}{n^{s}},\quad {\rm Re}\, s>0, s\neq 1.\label{2}
\end{align}
This leads to an integral representation of the zeta function in terms of the Mellin transform of the Fermi-Dirac distribution function \cite{abramowitz}. In this way, Riemann had defined the zeta function on the right half of the complex parameter space except for the point $s=1$. In order to extend it to the left half of the plane, Riemann derived two equivalent functional relations (we give only one, the other can be obtained from this by letting $s\rightarrow 1-s$)
\begin{equation}
\zeta(s) = 2 (2\pi)^{s-1}\sin \frac{\pi s}{2}\,\Gamma(1-s)\zeta(1-s),\quad {\rm Re}\, s<1.\label{3}
\end{equation}
This, therefore, generalizes the zeta function to the entire complex plane.

The zeta function is an analytic function \cite{edwards, titchmarsh} which has a simple pole at $s=1$ with residue $1$ (the pole structure is already manifest in \eqref{2}). It vanishes for $s=-2k$ where $k\geq 1$ which can be seen from \eqref{3}. These zeros are known as ``trivial" zeros of the zeta function since they arise from the kinematic trigonometric factor in \eqref{3} and are the only zeros for ${\rm Re}\, s \leq 0$ ($\zeta(0)=-\frac{1}{2}$). Furthermore, from the representation of the zeta function in terms of products involving prime numbers, it can be shown that the zeta function has no zero for ${\rm Re}\, s>1$ (since none of the factors can vanish there). Therefore, any other ``non trivial" zero of the zeta function must lie in the strip 
\begin{equation}
0<{\rm Re}\,s<1\label{4}.
\end{equation} 
Riemann conjectured \cite{riemann} that all other zeros of the zeta function lie on the critical line ${\rm Re}\,s=\frac{1}{2}$, namely,
\begin{equation}
\zeta ({\sc \frac{1}{2}} + i\lambda_{*}) = 0,\label{5}
\end{equation}
where $\lambda_{*}$ denotes the location of a zero on the critical line. This is known as the Riemann hypothesis and so far many zeros have been calculated on the critical line numerically \cite{odlyzko1, odlyzko2}. Note that the Riemann hypothesis specifies only the real part of the parameter $s$ at a zero. It does not say anything about the location (imaginary part of the parameter $s$) on the critical line. The numerical calculations do not yet show any particular recurrence relation or regularity in the locations of the zeros.

There is a fascinatingly beautiful symmetry in physics known as supersymmetry \cite {susy1,susy2, susy3,susy4}. It has a simple manifestation in one dimensional supersymmetric quantum mechanics \cite{susyqm1, susyqm2,susyqm3,susyqm4} where there are two conserved charges, $Q, Q^{\dagger}$ satisfying the graded algebra (the bracket with a subscript $+$ denotes an anti-commutator)
\begin{equation}
[Q, Q^{\dagger}]_{+} = H,\quad [Q,H] = 0 = [Q^{\dagger}, H],\label{6}
\end{equation}
where $H$ denotes the Hamiltonian of the quantum mechanical system and the adjoint corresponds to the Dirac adjoint (in our entire discussion). The Hamiltonian is assumed to be self-adjoint for unitarity of time evolution. An immediate consequence of the symmetry algebra in \eqref{6} is that the ground state energy of the theory vanishes, $E_{0}=0$ and all the other eigenstates of the Hamiltonian have real and positive energy,
\begin{equation}
E_{n} = \frac{\langle \psi_{n}|H|\psi_{n}\rangle}{\langle \psi_{n}|\psi_{n}\rangle} = \frac{(|Q |\psi_{n}\rangle|^{2}+|Q^{\dagger}|\psi_{n}\rangle|^{2})}{\langle\psi_{n}|\psi_{n}\rangle} \geq 0.\label{6a}
\end{equation}
A simple representation for the algebra \eqref{6} can be given in terms of $2\times 2$ matrices as
\begin{align}
Q & = \begin{pmatrix}
0 & 0\\
A & 0
\end{pmatrix},\quad Q^{\dagger} = \begin{pmatrix}
0 & A^{\dagger}\\
0 & 0
\end{pmatrix},\notag\\ 
H & = \begin{pmatrix}
H_{-} & 0\\
0 & H_{+}
\end{pmatrix} = \begin{pmatrix}
A^{\dagger} A & 0\\
0 & AA^{\dagger}
\end{pmatrix},\label{7}
\end{align}
where $A, A^{\dagger}$ denote the lowering and raising operators for the system. The Hamiltonians 
\begin{equation}
H_{-} = A^{\dagger}A,\quad H_{+}=AA^{\dagger},\label{8} 
\end{equation}
are known as supersymmetric partner Hamiltonians. The two Hamiltonians are almost isospectral in the sense that they share all the energy eigenvalues except for the ground state energy of $H_{-}$ which vanishes. Namely, if $|\psi_{0}\rangle$ denotes the ground state of $H_{-}$,
\begin{equation}
A|\psi_{0}\rangle = 0,\quad H_{-}|\psi_{0}\rangle = A^{\dagger}A |\psi_{0}\rangle = 0.\label{9}
\end{equation}
Furthermore, the supersymmetric partner states satisfy
\begin{align}
& |\psi_{n}\rangle = (A^{\dagger})^{n}|\psi_{0}\rangle,\quad  |\tilde{\psi}_{n}\rangle = A|\psi_{n}\rangle,\quad |\tilde{\psi}_{0}\rangle=0,\notag\\
& H_{-}|\psi_{n}\rangle = E_{n}|\psi_{n}\rangle,\ H_{+}|\tilde{\psi}_{n}\rangle = E_{n}|\tilde{\psi}_{n}\rangle.\label{9a}
\end{align}

There seems to be a tantalizing connection here. We are interested in studying the zeros of the zeta function on the critical line and the ground state energy of a supersymmetric system naturally vanishes. Therefore, if we can find a supersymmetric system whose energy eigenvalues are related to the zeta function, we can possibly understand the the zeros of the zeta function on the critical line ${\rm Re}\, s=\frac{1}{2}$ by looking at the ground state of the supersymmetric system. This is the idea which we would like to propose and pursue in this work to derive \eqref{5}.

We note that if an operator takes a monomial $x^{-s}$ to
\begin{equation}
\sum_{n=1}^{\infty} (-1)^{n+1} (nx)^{-s} = (1-2^{1-s})\zeta(s) x^{-s},\ {\rm Re}\,s>0,\label{10}
\end{equation}
where we have used \eqref{2}, we can obtain the zeta function as an eigenvalue of an operator acting on a monomial $x^{-s}, {\rm Re}\,s>0$. This basically involves a finite scale transformation. We recall that $x\frac{d}{dx}$, the generator of infinitesimal scale transformations, simply counts the power of $x$ in a monomial, namely, 
\begin{equation}
x\frac{dx^{-s}}{dx} = -sx^{-s},\quad (f(x\frac{d}{dx})x^{-s}) = f(-s)x^{-s}.\label{10a}
\end{equation} 
Therefore, we consider the operators 
\begin{align}
O & = \sum_{n=1}^{\infty} (-1)^{n+1} \exp\left((\ln n)\, x \frac{d}{dx}\right),\notag\\
O^{\dagger} & = \sum_{n=1}^{\infty}\frac{(-1)^{n+1}}{n} \exp\left((\ln n^{-1})\, x\frac{d}{dx}\right),\label{11}
\end{align}
where $O^{\dagger}$ is the Dirac adjoint of $O$. Using \eqref{10a}, it follows now that (see \eqref{2})
\begin{align}
O x^{-s} & = \sum_{n=1}^{\infty} (-1)^{n+1} e^{-s\ln n} x^{-s} = \sum_{n=1}^{\infty} \frac{(-1)^{n+1}}{n^{s}} x^{-s}\notag\\
& = (1-2^{1-s})\zeta(s) x^{-s},\ {\rm Re}\,s>0,\notag\\ 
O^{\dagger} x^{-s} & = \sum_{n=1}^{\infty} \frac{(-1)^{n+1}}{n} e^{s\ln n} x^{-s} = \sum_{n=1}^{\infty} \frac{(-1)^{n+1}}{n^{1-s}} x^{-s}\notag\\
&= (1-2^{s}) \zeta(1-s) x^{-s},\ {\rm Re}\, s < 1 .\label{12}
\end{align}
The action of the operators $O, O^{\dagger}$, therefore, naturally restricts the (negative) power of $x$ to lie in the strip \eqref{4} if we want to obtain zeta functions (as eigenvalues) and, in this strip, the operators also act on $|x|^{-s}$ exactly as in \eqref{12}. (The reason for choosing $|x|^{-s}$ will become clear when we discuss the normalizability of these functions.) We choose to work with wavefunctions in the strip of the general form $|x|^{-s} = |x|^{-\sigma+i\rho},\, 0<\sigma<1$ where $\sigma,\rho$ are real parameters.

The operators $O$ and $O^{\dagger}$ defined in \eqref{11} commute since the generators of scale transformations $x\frac{d}{dx}$ do. Therefore, we define the basic lowering and raising operators of the system as
\begin{equation}
A(\omega) = |x|^{-\frac{i\omega}{2}} O\, |x|^{-\frac{i\omega}{2}},\quad A^{\dagger}(\omega) = |x|^{\frac{i\omega}{2}} O^{\dagger}\, |x|^{\frac{i\omega}{2}},\label{13}
\end{equation}
where $\omega$ is a real parameter. We do not allow the parameter $\omega$ to be complex so that acting on a wavefunction in the strip, it does not change the domain of the wavefunction in the parameter space. More importantly, a complex $\omega$ would lead to a Hamiltonian which is not self-adjoint as we will discuss later. As a result, we can construct two supersymmetric partner Hamiltonians as (see \eqref{8})
\begin{align}
H_{-} & = A^{\dagger}(\omega)A(\omega) = |x|^{\frac{i\omega}{2}}O^{\dagger}O\,|x|^{-\frac{i\omega}{2}},\notag\\
H_{+} & = A(\omega) A^{\dagger}(\omega) = |x|^{-\frac{i\omega}{2}} OO^{\dagger}\,|x|^{\frac{i\omega}{2}}.\label{14}
\end{align}
Since this pair of Hamiltonians (formally) define a supersymmetric system, as discussed in \eqref{9}, the ground state energy of $H_{-}$ has to be zero. (We point out here that the infinitesmial scale generator $x\frac{d}{dx}$ and a variant of this have  been used as  Hamiltonians \cite{berry1,berry2,bender} in earlier studies of Riemann zeros from different perspectives, see also \cite{brody}. Here the basic element in our construction of $H$ is the group of finite scale transformations.)

Using \eqref{2} and \eqref{12}, it can be checked that acting on the space of functions $|x|^{-\sigma+i\rho}$, the lowering and raising operators give
\begin{align}
A(\omega) |x|^{-\sigma + i\rho} & = (1-2^{(1-\sigma)+i(\rho-{\sc \frac{\omega}{2}})})\zeta(\sigma-i(\rho-{\sc \frac{\omega}{2}}))\notag\\
&\qquad\times |x|^{-\sigma + i(\rho-\omega)},\notag\\
A^{\dagger}(\omega) |x|^{-\sigma + i\rho} & = (1-2^{\sigma-i(\rho+{\sc \frac{\omega}{2}})})\zeta(1-\sigma+i(\rho+{\sc \frac{\omega}{2}}))\notag\\
&\qquad \times|x|^{-\sigma + i(\rho+\omega)}.\label{15}
\end{align}
Namely, $A$ and $A^{\dagger}$ translate only the imaginary part of the exponent in the function by an amount $\mp \omega$ respectively up to overall multiplicative constant factors. As a result, we can think of $\sigma$ as a fixed constant for a given class of functions (namely, for a given $\sigma$, we can think of the functions as defined on a vertical line in the complex parameter space in the strip $0<\sigma<1$). Equation \eqref{15} leads to
\begin{widetext}
\begin{align}
H_{-} |x|^{-\sigma + i\rho} & = (1-2^{\sigma-i(\rho-{\sc \frac{\omega}{2}})})(1-2^{1-\sigma+i(\rho-{\sc \frac{\omega}{2}})})\zeta(\sigma-i(\rho-{\sc\frac{\omega}{2}})) \zeta(1-\sigma+i(\rho-{\sc\frac{\omega}{2}}))\, |x|^{-\sigma + i\rho},\notag\\
H_{+} |x|^{-\sigma + i\rho} & = (1-2^{\sigma-i(\rho+{\sc \frac{\omega}{2}})})(1-2^{1-\sigma+i(\rho+{\sc \frac{\omega}{2}})})\zeta(\sigma-i(\rho+{\sc\frac{\omega}{2}})) \zeta(1-\sigma+i(\rho+{\sc\frac{\omega}{2}}))\, |x|^{-\sigma + i\rho}.\label{16}
\end{align}
\end{widetext}
This shows that this class of functions define eigenfunctions of the partner Hamiltonians with products of zeta functions (and other factors) as eigenvalues. This is indeed what we started out looking for. 

The zeta functions in \eqref{15} and \eqref{16} are defined on the strip \eqref{4} in the complex parameter space, namely, $0<\sigma<1$. Therefore, it would seem that the ground state and the zero of the zeta function (see \eqref{9}) can, in general, lie anywhere in this strip. However, we have not yet discussed either the normalizability of the class of functions or the self-adjoint properties of the Hamiltonians in the Dirac sense  which we do next.  It is well known \cite{okikiolu, robertson,schaefer,bourbaki,endres, siera1} that the class of functions $x^{-\sigma+i\rho}$ is normalizable in the positive real axis $\mathbbm{R}_{+}$ in the Dirac sense (corresponding to a Dirac inner product) only for $\sigma=\frac{1}{2}$. Namely, only for the class of functions $\psi_{\rho}(x) = x^{-\frac{1}{2}+i\rho}$, we have
\begin{align}
\int_{0}^{\infty} dx\,\psi_{\rho}^{*}(x) \psi_{\rho'}(x) & = \int_{0}^{\infty} dx\, x^{-1-i(\rho-\rho')}\notag\\
& = 2\pi \delta(\rho-\rho').\label{17}
\end{align}
This particular class of functions is also complete in the Dirac sense, namely, for $x,y>0$,
\begin{align}
\int_{-\infty}^{\infty} d\rho\, \psi_{\rho}(x) \psi^{*}_{\rho} (y) & = \int_{-\infty}^{\infty} d\rho\, x^{-\frac{1}{2} + i\rho}\, y^{-\frac{1}{2} - i\rho}\notag\\
& = 2\pi \delta(x-y).\label{18}
\end{align}
Furthermore, these functions can be naturally extended to the negative axis leading to two linearly independent functions on the entire real axis which are normalizable and complete (in the Dirac sense)
\begin{align}
\psi_{even, \rho}(x) & = \frac{1}{\sqrt{4\pi}}\,|x|^{-\frac{1}{2}+i\rho},\notag\\
\psi_{odd, \rho} (x) & = \frac{{\rm sgn} (x)}{\sqrt{4\pi}}\,|x|^{-\frac{1}{2} + i\rho}.\label{19}
\end{align}
On the other hand, since there is no operator in our theory which can change the parity of a function, we can restrict ourselves to only one of the two classes of functions given in \eqref{19}. For simplicity (and because we are looking at the ground state) we choose to work with only the even class of functions in \eqref{19} (the other choice would also work equally well). Therefore, of the general class of functions we have considered in \eqref{15} and \eqref{16}, the normalizable functions correspond only to the choice $\sigma=\frac{1}{2}$. 

Next let us analyze whether the partner Hamiltonians are self-adjoint on the general class of functions $|x|^{-\sigma+i\rho}, 0<\sigma<1$. This is essential both for unitary time evolution as well as for supersymmetry in the system. Using \eqref{16} we can calculate
\begin{widetext}
\begin{align}
& \int dx\, (H_{-}\psi_{\rho}(x))^{*} \psi_{\rho'}(x) = (1-2^{\sigma + i(\rho-\frac{\omega}{2})} (1-2^{1-\sigma-i(\rho - \frac{\omega}{2}))}
\zeta (\sigma+i(\rho-\frac{\omega}{2}))\zeta (1-\sigma-i(\rho-\frac{\omega}{2}))\int dx\,|x|^{-2\sigma-i(\rho-\rho')},\notag\\
& \int dx\, \psi_{\rho}^{*}(x) (H_{-}\psi_{\rho'}(x)) = (1-2^{\sigma - i(\rho'-\frac{\omega}{2})} (1-2^{1-\sigma+i(\rho' - \frac{\omega}{2}))}\zeta (\sigma-i(\rho'-\frac{\omega}{2}))\zeta (1-\sigma+i(\rho'-\frac{\omega}{2}))\int dx\,|x|^{-2\sigma-i(\rho-\rho')}.\label{19a}
\end{align}
\end{widetext}
The two expressions are not the same even if we set $\rho'=\rho$ unless $\sigma=\frac{1}{2}$. Therefore, the Hamiltonian $H_{-}$ is not self-adjoint on a general class of functions $|x|^{-\sigma+i\rho}$ unless $\sigma = \frac{1}{2}$. We can come to the same conclusion for $H_{+}$ as well. This is already reflected in the energy eigenvalues in \eqref{16} being complex for a general $0<\sigma<1$ and can not be written as an absolute square, as supersymmetry would require (see \eqref{6a}), unless $\sigma=\frac{1}{2}$. This analysis shows that unitary time evolution as well as supersymmetry can not be realized on the general class of functions on the strip unless $\sigma= \frac{1}{2}$. (Incidentally, this analysis can also be carried out for a complex $\omega$ which would show that the partner Hamiltonians will be self-adjoint only if $\omega^{*} = \omega$ and $\sigma=\frac{1}{2}$.  Otherwise, unitary time evolution and supersymmetry will be violated. This is the main reason for choosing $\omega$ to be real.)

Normalizability, unitary time evolution as well as supersymmetry, therefore, select the space of functions to be of the form $|x|^{-\frac{1}{2}+i\rho}$ in the strip so that, for given values of $\rho$ and $\omega$, the ground state of the supersymmetric system has to satisfy (see \eqref{9})
\begin{widetext}
\begin{align}
A(\omega)|x|^{-\frac{1}{2}+i\rho} & = (1- 2^{\frac{1}{2}+i(\rho-\frac{\omega}{2})})\zeta(\frac{1}{2} - i(\rho-\frac{\omega}{2})) |x|^{-\frac{1}{2}+i(\rho-\omega)} = 0,\notag\\
H_{-} |x|^{-\frac{1}{2}+i\rho} & = |1-2^{\frac{1}{2}-i(\rho-\frac{\omega}{2})}|^{2} |\zeta(\frac{1}{2}-i(\rho-\frac{\omega}{2}))|^{2} |x|^{-\frac{1}{2}+i\rho} = 0.\label{20}
\end{align}
\end{widetext}
Here we have used the fact that $(\zeta(\frac{1}{2}-ix))^{*} = \zeta(\frac{1}{2}+ix)$. For \eqref{20} to be true, we must have 
\begin{equation}
\zeta (\frac{1}{2} -i(\rho-\frac{\omega}{2})) = \zeta(\frac{1}{2}+i\lambda_{*}) = 0,\label{21}
\end{equation}
where $\frac{\omega}{2}-\rho = \lambda_{*}$ corresponds to the location of a zero on the critical line, \eqref{5}, here arising naturally from the vanishing of ground state energy in a supersymmetric quantum mechanical model (within the standard quantum mechanical framework). We note that if $\frac{\omega}{2}-\rho=\lambda$ does not coincide with the location of a zero of the zeta function on the critical line, then $|\zeta(\frac{1}{2}+i\lambda)| \neq 0$ and the wavefunction would correspond to a positive energy state as supersymmetry would require. Namely, being an absolute square, the energy eigenvalues in \eqref{20} will always be positive (even when the zeta function is negative) unless it is zero. We can now write explicitly the ground state wavefunction of $H_{-}$ (not normalized) to be
\begin{equation}
\psi_{0}(x) = |x|^{-\frac{1}{2} + i(\frac{\omega}{2}-\lambda_{*})},\ A\psi_{0}(x)=0,\ E_{0}=0.\label{22}
\end{equation}
We can build the Hilbert space of the supersymmetric theory on this ground state using \eqref{9a}, leading to
\begin{align}
\psi_{n}(x) & = C_{n} |x|^{-\frac{1}{2} + i(n\omega +(\frac{\omega}{2} - \lambda_{*}))},\ n=1,2, \cdots,\notag\\
\tilde{\psi}_{n}(x) & = \tilde{C}_{n} |x|^{-\frac{1}{2} + i(n\omega - (\frac{\omega}{2}+\lambda_{*}))},\quad \tilde{\psi}_{0}(x)=0,\label{23}
\end{align}
where the constants $C_{n}$ and $\tilde{C}_{n}$ are given by
\begin{align}
C_{n} & = \prod_{m=1}^{n} (1- 2^{\frac{1}{2} -i(m\omega-\lambda_{*})}) \zeta (\frac{1}{2}+i(m\omega-\lambda_{*})),\notag\\
\tilde{C}_{n} & = |(1-2^{\frac{1}{2}-i(n\omega -\lambda_{*})}) \zeta(\frac{1}{2}+i(\lambda_{*}-n\omega))|^{2} C_{n-1}.\label{24}
\end{align}
Both the partner states, $\psi_{n}(x), \tilde{\psi}_{n}(x)$, share the same energy eigenvalue
\begin{equation}
E_{n} = |(1-2^{\frac{1}{2}-i(n\omega-\lambda_{*})})|^{2}|\zeta(\frac{1}{2}+i(\lambda_{*}-n\omega))|^{2},\label{25}
\end{equation}
for $n=1,2,\cdots$, as supersymmetry will require. We also note from \eqref{22} that since $x\frac{d}{dx}$ counts the power of $x$ in a monomial (see \eqref{10a}),
\begin{equation}
x\frac{d\psi_{0}(x)}{dx} = (-\frac{1}{2} + i(\frac{\omega}{2} -\lambda_{*}))\psi_{0}(x).\label{25a}
\end{equation}
Such an equation was already proposed in \cite{berry1} from a different perspective. In our case, through a similarity transformation, this can actually be rewritten as an eigenvalue equation for the ground state with the eigenvalue given by the location of the zero as
\begin{equation}
B\psi_{0}(x) = \lambda_{*} \psi_{0}(x),\ \ B = i |x|^{-\frac{(1-i\omega)}{2}} x\frac{d}{dx} |x|^{\frac{(1-i\omega)}{2}}.\label{25b}
\end{equation}

The energy eigenstates in \eqref{22} and \eqref{23} define a discrete set of basis states of the form $|x|^{-\frac{1}{2}+i(C+n\omega)}$ where we have identified $C= \pm\frac{\omega}{2}-\lambda_{*}$. Therefore, we briefly indicate how normalizability and completeness hold for these states. Identifying these basis states on the positive real axis as $\phi_{n} (x), n=0,1,2,\cdots$ we note that
\begin{equation}
\int dx\,\phi_{n}^{*}(x)\phi_{n'}(x) = \int dx\, x^{-1- i(n-n')\omega}.\label{26}
\end{equation} 
If we now define a complex variable $z=x^{i\omega}$, $\phi_{n}(x)\rightarrow\phi_{n}(z)= z^{\frac{1}{\omega}(\frac{i}{2}+C+n\omega)}$  and the integral in \eqref{26} leads to
\begin{equation}
\int dx\,\phi_{n}^{*}(x)\phi_{n'}(x) = -\frac{i}{\omega}\oint \frac{dz}{z^{1+(n-n')}}=\frac{2\pi}{\omega}\delta_{nn'},\label{27}
\end{equation}
where the contour integral is taken along a unit circle around the origin in an anti-clockwise direction. The functions
\begin{equation}
\phi_{n}(z) = z^{\frac{1}{\omega} (\frac{i}{2}+C)} z^{n},\label{28}
\end{equation}
define a complete basis for any function of the form $z^{\frac{1}{\omega} (\frac{i}{2}+C)} f(z)$ where $f(z)$ is analytic at the origin, which can be seen simply from
\begin{equation}
z^{\frac{1}{\omega} (\frac{i}{2}+C)} f(z) = \sum_{n=0}^{\infty} c_{n}\phi_{n}(z),\label{29}
\end{equation}
where $c_{n} = \frac{f^{(n)}(0)}{n!}$.

We note that the Riemann zeta function has an infinite number of ``nontrivial" zeros on the critical line which may seem to indicate that our theory has an infinite degeneracy of ground states. In reality, however, the parameter $\omega$ which defines our theory can be carefully chosen so that the theory indeed has a unique vacuum state. For different values of $\omega$ (different theories) one can realize the other Riemann zeros as the  ground state energy of the corresponding theories.

In summary, we have proposed a new way of looking at the Riemann zeros on the critical line by looking at the vanishing ground state energy in a supersymmetric theory (within standard quantum framework). In a simple model, this naturally leads to \eqref{5}. In the following appendix, we give some more details of some of the essential derivations.

We would like to thank Drs  Ashok Kumar Diktiya and Levi Greenwood for discussions.

\appendix

\section{Details of some of the derivations}

The derivations given in the paper are quite  self-complete. However, here we give more details of some of the essential equations in the paper in order to help the readers. 

\begin{enumerate}

\item Using  \eqref{10a} and \eqref{11} in the paper,  \eqref{12} follows in a straightforward manner. Here is a detailed derivation of it. We note that

\begin{align}
O x^{-s} & =  \sum_{n=1}^{\infty}(-1)^{n+1}\exp((\ln n)x\frac{d}{dx})\, x^{-s}\nonumber \\
& = \sum_{n=1}^{\infty}(-1)^{n+1}\sum_{m=0}^{\infty}\frac{1}{m!} (\ln n)^{m} \left(x\frac{d}{dx}\right)^{m} x^{-s}\notag\\
& = \sum_{n=1}^{\infty} (-1)^{n+1} \sum_{m=0}^{\infty} \frac{(\ln n)^{m} (-s)^{m}}{m!}\notag\\
& = \sum_{n=1}^{\infty} (-1)^{n+1} \exp\left(-s\ln n\right) x^{-s} = \sum_{n=1}^{\infty}\frac{(-1)^{n+1}}{n^{s}} x^{-s}\notag\\
& = (1-2^{1-s})\zeta(s) x^{-s},\quad {\rm Re}\, s>0.\label{app1}
\end{align}
where in the last line, we have used the definition of the Riemann zeta function given in  \eqref{2} of the paper. In the same manner, it can be shown that
\begin{align}
O^{\dagger} x^{-s} & = \sum_{n=1}^{\infty}\frac{(-1)^{n+1}}{n}\exp((\ln n^{-1})x\frac{d}{dx}) x^{-s}\nonumber \\
& = \sum_{n=1}^{\infty}\frac{(-1)^{n+1}}{n}\sum_{m=0}^{\infty}\frac{1}{m!} (\ln n^{-1})^{m} \left(x\frac{d}{dx}\right)^{m} x^{-s}\notag\\
& = \sum_{n=1}^{\infty} \frac{(-1)^{n+1}}{n} \sum_{m=0}^{\infty} \frac{(\ln n^{-1})^{m} (-s)^{m}}{m!} x^{-s}\notag\\
& = \sum_{n=1}^{\infty} \frac{(-1)^{n+1}}{n} \exp\left(-s\ln n^{-1}\right) x^{-s}\notag\\
& = \sum_{n=1}^{\infty}\frac{(-1)^{n+1}}{n^{1-s}} x^{-s}\notag\\
& = (1-2^{s})\zeta(1-s) x^{-s},\quad {\rm Re}\, s <1.\label{app2}
\end{align}

\item Derivation of Eq. \eqref{15} in the paper: Recall that \eqref{app1} and \eqref{app2} above hold also for $|x|^{-s} = |x|^{-\sigma + i\rho}$, as mentioned in the paper, in the strip $0< {\rm Re}\, s = \sigma <1$. Therefore, using \eqref{app1} and \eqref{app2} above, we obtain

\begin{align}
\lefteqn{A(\omega)|x|^{-\sigma+i\rho}= |x|^{-\frac{i\omega}{2}}O|x|^{-\frac{i\omega}{2}}|x|^{-\sigma+i\rho}}\nonumber \\
& = |x|^{-\frac{i\omega}{2}} O|x|^{-\sigma+i(\rho-\frac{\omega}{2})}\nonumber \\
& = (1-2^{1-\sigma+i(\rho-\frac{\omega}{2})})\zeta (\sigma-i(\rho-\frac{\omega}{2}))|x|^{-\sigma+i(\rho-\omega)},\label{app3}\\
\lefteqn{A^{\dagger}(\omega)|x|^{-\sigma+i\rho}=  |x|^{\frac{i\omega}{2}}O^{\dagger} |x|^{\frac{i\omega}{2}}|x|^{-\sigma+i\rho}}\nonumber \\
& = |x|^{\frac{i\omega}{2}}O^{\dagger}|x|^{-\sigma+i(\rho+\frac{\omega}{2})}\nonumber \\
& = (1-2^{\sigma-i(\rho+\frac{\omega}{2})})\zeta (1-\sigma+i(\rho+\frac{\omega}{2}))|x|^{-\sigma+i(\rho+\omega)}.\label{app4}
\end{align}
In deriving \eqref{app3} and \eqref{app4},  we have used our earlier results derived in \eqref{app1} and \eqref{app2}.

\item Derivation of \eqref{16} in the paper: This basically follows from \eqref{app3} and \eqref{app4} derived above. Note that $A(\omega)$ lowers the imaginary part of the power of $|x|$ by $\omega$ while $A^{\dagger}(\omega)$ raises the imaginary part of the power by $\omega$ (besides giving some multiplicative constants in the function)

We now see how $H_{-}$ operates on $|x|^{-\sigma+i\rho}$, 
\begin{align}
\lefteqn{H_{-}|x|^{-\sigma+i\rho}  = A^{\dagger}(\omega)A(\omega) |x|^{-\sigma+i\rho}}\notag\\
& = (1-2^{1-\sigma+i(\rho-\frac{\omega}{2})})\zeta (\sigma-i(\rho-\frac{\omega}{2})) \notag\\
& \qquad \times A^{\dagger}(\omega)|x|^{-\sigma+i(\rho-\omega)}\notag\\
& = (1-2^{\sigma-i(\rho -\frac{\omega}{2})})(1-2^{1-\sigma+i(\rho-\frac{\omega}{2}})\notag\\
&\qquad \times \zeta(\sigma - i(\rho-\frac{\omega}{2}))\zeta(1-\sigma + i(\rho-\frac{\omega}{2}))|x|^{-\sigma+i\rho}.\label{app5}
\end{align}
In exactly the same manner, we obtain
\begin{align}
\lefteqn{H_{+}|x|^{-\sigma+i\rho} = A(\omega)A^{\dagger}(\omega) |x|^{-\sigma+i\rho}}\notag\\
& = (1-2^{\sigma-i(\rho+\frac{\omega}{2})})\zeta (1-\sigma+i(\rho+\frac{\omega}{2}))\notag\\
&\qquad \times A(\omega) |x|^{-\sigma+i(\rho+\omega)}\notag\\
& = (1-2^{\sigma-i\left(\rho+\frac{\omega}{2}\right)})(1-2^{1-\sigma+i\left(\rho+\frac{\omega}{2}\right)})\notag\\
&\qquad\times \zeta(\sigma-i(\rho+\frac{\omega}{2}))\zeta(1-\sigma+i(\rho+\frac{\omega}{2}))|x|^{-\sigma+i\rho}.\label{app6}
\end{align}

\item Explaining \eqref{19a} in the paper: A Hamiltonian (or any other operator) is self-adjoint on a space of functions if it satisfies the condition

$$
\langle Hg|f\rangle = \langle g|Hf\rangle,
$$
where $|f\rangle, |g\rangle$ are two states in the Hilbert space. In our case, in \eqref{19a}, the states $|\psi_{\rho}\rangle , |\psi_{\rho'}\rangle$ correspond to the state functions

\begin{align*}
\psi_{\rho}(x) & = \langle x|\psi_{\rho}\rangle = |x|^{-\sigma+i\rho},\\
\psi_{\rho'}(x) & = \langle x|\psi_{\rho'}\rangle = |x|^{-\sigma+i\rho'}.
\end{align*}
Therefore, Eq. (23) which actually follows from \eqref{app5} above, simply says that

$$
\langle H_{-}\psi_{\rho}|\psi_{\rho'}\rangle \neq \langle\psi_{\rho}|H_{-}\psi_{\rho'}\rangle.
$$
unless $\sigma=\frac{1}{2}$ and, consequently, it is not self-adjoint otherwise. This is reflected in the fact that the energy eigenvalues in \eqref{app5} are complex if $\sigma\neq \frac{1}{2}$. The same conclusion follows for $H_{+}$ as well.

\item Higher states of the system: As we have emphasized $A^{\dagger}(\omega)$ acting on a function $|x|^{-\sigma+i\rho}$ raises the imaginary part of the exponent of $|x|$ by $\omega$ and gives a multiplicative factor as given in \eqref{app4} above. Since the ground state of $H_{-}$ has the form $|x|^{-\frac{1}{2} +i(\frac{\omega}{2}-\lambda_{*})}$, by repeated application of the raising operator $n$ times, we obtain the $n$th higher state of $H_{-}$ to be

$$
\psi_{n}(x) = (A^{\dagger}(\omega))^{n} \psi_{0}(x) = C_{n} |x|^{-\frac{1}{2}+i(n\omega + (\frac{\omega}{2}-\lambda_{*}))}.
$$
The supersymmetric partner state (eigenstate of $H_{+}$) is obtained by applying the lowering operator

$$
\tilde{\psi}_{n} (x) = A(\omega) \psi_{n}(x).
$$
Since the lowering operator lowers the imaginary part of the exponent of $|x|$ by $\omega$ (and gives some multiplicative constant as shown in \eqref{app3} above), we obtain

$$
\tilde{\psi}_{n}(x) = \tilde{C}_{n} |x|^{-\frac{1}{2}+i(n\omega-(\frac{\omega}{2}+\lambda_{*}))},
$$
where $C_{n}, \tilde{C}_{n}$ are multiplicative constants given in \eqref{24} of the paper.

\end{enumerate}

\end{document}